\documentclass[11pt,a4paper]{article}
\title{\bf A unified treatment for ODEs under Osgood and Sobolev type conditions}
\author{Huaiqian \uppercase{LI},
\quad Dejun \uppercase{LUO}\footnote{Email: luodj@amss.ac.cn.
The second author is supported by the National Natural Science
Foundation of China (Grant No. 11101407), and the Key Laboratory of Random Complex Structures
and Data Science, Academy of Mathematics and Systems Science, Chinese Academy of Sciences (No. 2008DP173182).}
\vspace{3mm}\\
{\footnotesize Institute of Applied Mathematics, Academy of Mathematics and Systems Science,}\\
{\footnotesize Chinese Academy of Sciences, Beijing 100190, China} }
\date{}
\usepackage{amssymb,amsmath,amsfonts,amsthm,array,bm,color}

\def\R{\mathbb{R}}

\def\L{\mathcal{L}}

\newcommand{\ra}{\rightarrow}

\newcommand{\da}{\downarrow}

\newcommand{\ee}{\varepsilon}

\def\div{\textup{div}}
\def\d{\textup{d}}

\def\fin{\hfill$\square$}
\def\<{\langle}
\def\>{\rangle}

\def\Xint#1{\mathchoice
{\XXint\displaystyle\textstyle{#1}}%
{\XXint\textstyle\scriptstyle{#1}}%
{\XXint\scriptstyle\scriptscriptstyle{#1}}%
{\XXint\scriptscriptstyle\scriptscriptstyle{#1}}%
 \!\int}
\def\XXint#1#2#3{{\setbox0=\hbox{$#1{#2#3}{\int}$}
\vcenter{\hbox{$#2#3$}}\kern-.5\wd0}}

\def\dashint{\Xint-}
\def\bint{\dashint}

\newtheorem{theorem}{Theorem}[section]
\newtheorem{lemma}[theorem]{Lemma}       
\newtheorem{corollary}[theorem]{Corollary}
\newtheorem{proposition}[theorem]{Proposition}
\newtheorem{remark}[theorem]{Remark}
\newtheorem{example}[theorem]{Example}
\newtheorem{definition}[theorem]{Definition}

\begin{document}

\maketitle
\makeatletter 
\renewcommand\theequation{\thesection.\arabic{equation}}
\@addtoreset{equation}{section}
\makeatother 

\begin{abstract}
In this paper we present a unified treatment for the ordinary differential
equations under the Osgood and Sobolev type conditions, following Crippa and de
Lellis's direct method. More precisely,
we prove the existence, uniqueness and regularity of
the DiPerna--Lions flow generated by a vector field which is
``almost everywhere Osgood continuous''.
\end{abstract}



\section{Introduction}

In the seminal paper \cite{DiPernaLions89}, DiPerna and Lions
established the existence and uniqueness of the quasi-invariant flow
of measurable maps generated by Sobolev vector fields with bounded
divergence. Their method is quite indirect in the sense that they
first established the well-posedness of the corresponding transport
equation, from which they deduced the results on ODE. Their
methodology is now called the DiPerna--Lions theory and can be regarded
as a generalization of the classical method of characteristics. It
has subsequently been extended to the case of BV vector fields with bounded
divergence by Ambrosio \cite{Ambrosio04}, via the well-posedness of
the continuity equation (cf. \cite{Ambrosio08,AmbrosioCrippa} for a
detailed account of these results). This theory was generalized in
\cite{AmbrosioFigalli09, FangLuo10} to the infinite dimensional Wiener space.
Recently, Crippa and de Lellis \cite{CrippadeLellis}
obtained some a-priori estimates on the flow (called regular
Lagrangian flow there) directly in Lagrangian formulation, namely
without exploiting the connection of the ODE with the transport or
continuity equation. Applying this method, they can give a direct
construction of the flow (see the extension to the case of
stochastic differential equations in \cite{Zhang10,
FangLuoThalmaier}).

To introduce the setting of the present work, we recall the key
ingredient in Crippa--de Lellis's direct method, namely, a Sobolev
vector field $b\in W_{loc}^{1,p}(\R^d)\, (p\geq1)$ is ``almost
everywhere Lipschitz continuous'' (it holds even for BV vector
fields, see \cite[Lemma A.3]{CrippadeLellis}). More precisely, there
are a negligible subset $N\subset \R^d$ and a constant $C_d$
depending only on the dimension $d$, such that for all $x,y\notin N$
and $|x-y|\leq R$, one has
  \begin{equation}\label{ae-Lipschitz}
  |b(x)-b(y)|\leq C_d|x-y|\big(M_R|\nabla b|(x)+M_R|\nabla
  b|(y)\big),
  \end{equation}
where $M_Rf$ is the local maximal function of $f\in
L_{loc}^1(\R^d)$:
  $$M_Rf(x):=\sup_{0<r\leq R}\frac1{\L^d(B(x,r))}\int_{B(x,r)}|f(y)|\,\d y.$$
Here $\L^d$ is the Lebesgue measure on $\R^d$ and $B(x,r)$ is the
ball centered at $x$ in $\R^d$ with radius $r>0$. If $x=0$ is the origin, we
will simply write $B(r)$ instead of $B(0,r)$. For a proof of the
inequality \eqref{ae-Lipschitz}, see
\cite[Appendix]{FangLuoThalmaier}; one can also find in \cite{Hajlasz} a
discussion on Sobolev spaces on general metric space. Using
the inequality \eqref{ae-Lipschitz}, Crippa and de Lellis estimated
the following type of quantity
  \begin{equation}\label{log-functional}
  \int_{B(R)}\log\bigg(\frac{|X_t(x)-\tilde X_t(x)|}\delta +1\bigg)\,\d x
  \end{equation}
in terms of $R, \delta$ and the $L^p$-norms of $\nabla b, \nabla
\tilde b$ on some ball, where $X_t$ and $\tilde X_t$ are
respectively the flows associated to the Sobolev vector fields $b$
and $\tilde b$. In a recent paper \cite{LiuHuang}, the authors slightly
generalized the results in \cite{CrippadeLellis} by assuming that
$|\nabla b|\in (L\log \cdots \log L)_{loc}$. We refer the readers to \cite{BouchutCrippa} for
some further extensions.

On the other hand, the study of stochastic differential equations
with non-Lipschitz coefficients has attracted intensive attentions
in the past decade; see for instance \cite{Malliavin, FangZhang,
Zhang05a, Luo11}. In particular, Fang and Zhang considered in
\cite{FangZhang} the general Osgood condition:
  \begin{equation}\label{Osgood}
  |b(x)-b(y)|\leq C|x-y|r(|x-y|^2),\quad |x-y|\leq c_0,
  \end{equation}
where $r:(0,c_0]\ra[1,\infty)$ is a continuous function defined on a
neighborhood of $0$ and satisfies $\int_0^{c_0}\frac{\d
s}{sr(s)}=\infty$. Under this condition and assuming that the ODE
  $$\frac{\d X_t}{\d t}=b(X_t),\quad X_0=x$$
has no explosion, they proved that the solution $X_t$ is a flow of
homeomorphisms on $\R^d$ (see \cite[Theorem 2.7]{FangZhang}). If in
addition $r(s)=\log\frac1s$ and the generalized divergence of $b$ is
bounded, then it is shown in \cite[Theorem 1.8]{Fang09} that the
Lebesgue measure $\L^d$ is also quasi-invariant under the flow
$X_t$. In a recent paper \cite{Luo09}, the second named author
extended this result to the Stratonovich SDE with smooth
diffusion coefficients, using Kunita's formula for the
Radon--Nikodym derivative of the stochastic flow (see \cite[Lemma
4.3.1]{Kunita}).

Inspired by these two types of conditions \eqref{ae-Lipschitz} and
\eqref{Osgood}, we consider in this work the following assumption on
the time dependent measurable vector field $b:[0,T]\times\R^d\ra
\R^d$:
\begin{enumerate}
\item[(H)] there exist a nonnegative function $g\in L^1\big([0,T],L^1_{loc}(\R^d)\big)$ and
negligible subsets $N_t$, such that for all $t\in[0,T]$ and
$x,y\notin N_t$, one has
  \begin{equation}\label{Osgood-Sobolev}
  |b_t(x)-b_t(y)|\leq \big(g_t(x)+g_t(y)\big)\rho(|x-y|),
  \end{equation}
where $\rho\in C(\R_+,\R_+)$ is strictly increasing, $\rho(0)=0$ and
$\int_{0+}\frac{\d s}{\rho(s)}=\infty$.
\end{enumerate}
The typical examples of the function $\rho$ are
$\rho(s)=s,\,s\log\frac1s,\,s\big(\log\frac1s\big)
\big(\log\log\frac1s\big),\cdots$. Notice that the latter two
functions are only well defined on some small interval $(0,c_0]$,
but we can extend their domain of definition by piecing them
together with radials. In this paper we fix such a function $\rho$.
Similarly we may call a function satisfying \eqref{Osgood-Sobolev}
``almost everywhere Osgood continuous'' (see the end of Section 2 for an example
of this kind of functions). It is clear that if we take
$g_t=C_d M_R|\nabla b_t|$ and $\rho(s)=s$ for all $s\geq0$, then the
inequality \eqref{Osgood-Sobolev} is reduced to
\eqref{ae-Lipschitz}. On the other hand, if $g$ is essentially
bounded, then \eqref{Osgood-Sobolev} becomes the general Osgood
condition \eqref{Osgood}, except on the negligible set $N_t$ (we can
redefine $b_t$ on this null set to get a continuous vector field).
Therefore, this paper can be seen as a unified treatment of the two
different types of conditions. We would like to mention that
Professor L. Ambrosio once told Professor S. Fang and the second author
(by private communication) that there might be no unified framework for the
DiPerna--Lions theory and the Osgood type conditions. The
assumptions like (H) were considered in \cite{RocknerZhang10}, but
the function $\rho$ was always taken as $\rho(s)=s$ for all
$s\geq 0$.

The paper is organized as follows. In Section 2, we recall the definition of the
regular Lagrangian flow and some important facts. An example of a vector field $b$
satisfying the condition (H) is also given there. Then in section 3, we construct a
unique regular Lagrangian flow under the condition (H) and the boundedness
of the divergence of $b$, following the direct
method of Crippa and de Lellis. Finally in Section 4, we show a regularity property of the
flow, which is weaker than the approximate differentiability
discussed in \cite{CrippadeLellis} and prove a compactness
result on the flow. To avoid technical complexities, we assume that
the vector fields are bounded throughout this paper.

\section{Preparations and an example}

We first give the definition of the flow associated to a vector
field $b$ (also called regular Lagrangian flow in \cite{Ambrosio04,
CrippadeLellis}).

\begin{definition}[Regular Lagrangian flow]\label{sect-2-def}
Let $b\in L^1_{loc}([0,T]\times\R^d,\R^d)$. A map
$X:[0,T]\times\R^d\ra\R^d$ is called the regular Lagrangian flow
associated to the vector field $b$ if
\begin{enumerate}
\item[\rm(i)] for a.e. $x\in\R^d$, the function $t\ra X_t(x)$ is
absolutely continuous and satisfies
  $$X_t(x)=x+\int_0^t b_s(X_s(x))\,\d s,\quad\mbox{for all }t\in[0,T];$$
\item[\rm(ii)] there exists a constant $L>0$ independent of $t\in[0,T]$ such
that $(X_t)_\#\L^d\leq L\L^d$.
\end{enumerate}
\end{definition}

Recall that $\L^d$ is the Lebesgue measure on $\R^d$ and
$(X_t)_\#\L^d$ is the push-forward of $\L^d$ by the flow $X_t$. $L$
will be called the compressibility constant of the flow $X$.

Next we introduce some notations and results that will be used in
the subsequent sections. Denote by $\Gamma_T=C([0,T],\R^d)$, i.e.
the space of continuous paths in $\R^d$. For $\gamma\in\Gamma_T$, we
write $\|\gamma\|_{\infty,T}$ for its supremum norm. Let $\delta>0$,
we define an auxiliary function by (cf. \cite[(2.7)]{FangZhang})
  \begin{equation}\label{psi-delta}
  \psi_\delta(\xi)=\int_0^\xi \frac{\d s}{\rho(s)+\delta},\quad \xi\geq0.
  \end{equation}
Note that if $\rho(s)=s$ for all $s\geq0$, then
  $$\psi_\delta(\xi)=\int_0^\xi \frac{\d s}{s+\delta}=\log\bigg(\frac \xi\delta+1\bigg)$$
which is the functional used in \eqref{log-functional}. Here are
some properties of $\psi_\delta$.

\begin{lemma}\label{sect-2-lem.1}
\begin{enumerate}
\item[\rm(1)] $\lim_{\delta\da0} \psi_\delta(\xi)=+\infty$ for all
$\xi>0$;
\item[\rm(2)] for any $\delta>0$, the function $\psi_\delta$ is
concave.
\end{enumerate}
\end{lemma}

\noindent{\bf Proof.} Property (1) follows from the fact that
$\int_{0+}\frac{\d s}{\rho(s)}=\infty$. To prove (2), we notice that
$\psi'_\delta(s)=\frac{1}{\rho(s)+\delta}$. Since $s\mapsto \rho(s)$
is increasing, we see that the derivative $\psi'_\delta(s)$ is
monotone decreasing, and hence $\psi_\delta$ is concave. \fin

\medskip

The concavity of $\psi_\delta$ will play an important role in the
arguments of Section 4. Finally we give an inequality concerning the
local maximal function (see \cite[Lemma A.2]{CrippadeLellis}).

\begin{lemma}\label{sect-2-lem.2}
Let $R,\,\lambda$ and $\alpha$ be positive constants. Then there is
$C_d$ depending only on the dimension $d$, such that
  $$\L^d\{x\in B(R):M_\lambda f(x)>\alpha\}\leq \frac{C_d}\alpha
  \int_{B(R+\lambda)}|f(y)|\,\d y.$$
\end{lemma}

Before finishing this section, we give an example of a vector field $b:\R^d\ra\R^d$
which satisfies the condition (H), but at the same time satisfies neither \eqref{ae-Lipschitz}
nor \eqref{Osgood}. The basic idea is to consider the sum of a Sobolev vector field
and an Osgood continuous vector field.

\begin{example}\label{sect-2-exa}
{\rm For $t\in\R$, let
  $$V(t)=\sum_{k=1}^\infty \frac{|\sin kt|}{k^2}.$$
Then by \cite[(2.12)]{FangZhang}, we have
  $$|V(t)-V(s)|\leq C_1|t-s|\log\frac1{|t-s|},\quad\mbox{for } |t-s|\leq e^{-1},$$
where $C_1>0$ is a constant. Define the function $\rho:\R_+\ra\R_+$ as follows:
  $$\rho(s)=\begin{cases}
  s\log\frac1s, & 0\leq s\leq e^{-2};\\
  s+e^{-2},& s>e^{-2}.
  \end{cases}$$
Then $\rho\in C^1(\R_+)$ is strictly increasing and $\rho(s)\geq s\log\frac1s$ for all $s\in[0,e^{-1}]$, hence
  $$|V(t)-V(s)|\leq C_1\rho(|t-s|),\quad\mbox{for } |t-s|\leq e^{-1}.$$
It is clear that $0\leq V(t)\leq \frac{\pi^2}6$, thus for $|t-s|>e^{-1}$,
  $$|V(t)-V(s)|\leq \frac{\pi^2}3\leq \frac{\pi^2 e}3\rho(e^{-1})
  \leq \frac{\pi^2 e}3\rho(|t-s|).$$
Let $C_2:= C_1\vee\frac{\pi^2 e}3$; then
  \begin{equation}\label{sect-2-exa.1}
  |V(t)-V(s)|\leq C_2\rho(|t-s|),\quad\mbox{for all } t,s\in\R.
  \end{equation}
As mentioned in \cite[Remark 2.10]{FangZhang}, $V(t)\sim t\log\frac1t$ as
$t\da0$, hence it is not locally Lipschitz continuous.

Now we take a Sobolev vector field $b_1\in W^{1,p}(\R^d,\R^d)\,(p>1)$ which is not
continuous. For $x\in\R^d$,
define $b_2(x)=(V(x_1),\ldots,V(x_d))$. Then by \cite[THEOREM 1]{Hajlasz},
there exist $g\in L^p(\R^d)$ and a negligible set $N\subset\R^d$ such that
  \begin{equation*}
  |b_1(x)-b_1(y)|\leq |x-y|(g(x)+g(y)),\quad\mbox{for all } x,y\notin N.
  \end{equation*}
Next by \eqref{sect-2-exa.1}, it holds for all $x,y\in\R^d$ that
  \begin{equation*}
  |b_2(x)-b_2(y)|\leq \sum_{i=1}^d |V(x_i)-V(y_i)|\leq \sum_{i=1}^d C_2 \rho(|x_i-y_i|)
  \leq C_2d\,\rho(|x-y|).
  \end{equation*}
Finally we set $b=b_1+b_2$. Then by the above two inequalities, for all $x,y\notin N$,
  \begin{align*}
  |b(x)-b(y)|&\leq |b_1(x)-b_1(y)|+|b_2(x)-b_2(y)|\cr
  &\leq |x-y|(g(x)+g(y))+C_2d\,\rho(|x-y|)\cr
  &\leq C_2d\,\rho(|x-y|)\big[(1+g(x))+(1+g(y))\big],
  \end{align*}
where in the last inequality we have used the facts that $C_2d>1$ and $\rho(s)\geq s$
for all $s\geq0$. Note that the function $1+g\in L^p_{loc}(\R^d)$, thus it is locally integrable.}
\end{example}

\section{Existence and uniqueness of the regular Lagrangian flow}

In order to prove the existence and uniqueness of the flow generated
by a vector field $b$ satisfying the assumption (H), we first
establish an a-priori estimate.

\begin{theorem}\label{a-priori estimate}
Let $b$ and $\tilde b$ be time dependent bounded vector fields
satisfying {\rm(H)} with $g$ and $\tilde g$ respectively. Let $X$
and $\tilde X$ be the regular Lagrangian flows associated to $b$ and
$\tilde b$, with the compressibility constants $L$ and $\tilde L$
respectively. Then for any $R>0$ and $t\leq T$,
  $$\int_{B(R)}\psi_\delta\big(\|X_\cdot(x)-\tilde X_\cdot(x)\|_{\infty,T}\big)\,\d x
  \leq (L+\tilde L)\|g\|_{L^1([0,T]\times B(\bar R))}
  +\frac{\tilde L}\delta \|b-\tilde b\|_{L^1([0,T]\times B(\bar R))},$$
where $\psi_\delta$ is defined in \eqref{psi-delta},
$\|\cdot\|_{\infty,T}$ is the supremum norm in $\Gamma_T$ and $\bar
R=R+T\big(\|b\|_{L^\infty}\vee\|\tilde b\|_{L^\infty}\big)$.
\end{theorem}

\noindent{\bf Proof.} By Definition \ref{sect-2-def}(1), for a.e.
$x\in\R^d$, the function $t\mapsto |X_t(x)-\tilde X_t(x)|$ is
Lipschitz continuous, hence
  \begin{align*}
  \frac{\d}{\d t}\big[\psi_\delta\big(|X_t(x)-\tilde X_t(x)|\big)\big]
  &=\psi'_\delta\big(|X_t(x)-\tilde X_t(x)|\big)
  \frac{\d}{\d t}|X_t(x)-\tilde X_t(x)|\cr
  &\leq\frac{\big|b_t(X_t(x))-\tilde b_t(\tilde X_t(x))\big|}
  {\rho\big(|X_t(x)-\tilde X_t(x)|\big)+\delta}.
  \end{align*}
Integrating from 0 to $t$ and noticing that $\psi_\delta(0)=0$, we
get
  $$\psi_\delta\big(|X_t(x)-\tilde X_t(x)|\big)
  \leq \int_0^t\frac{\big|b_s(X_s(x))-\tilde b_s(\tilde X_s(x))\big|}
  {\rho\big(|X_s(x)-\tilde X_s(x)|\big)+\delta}\,\d s.$$
As a result,
  $$\sup_{0\leq t\leq T}\psi_\delta\big(|X_t(x)-\tilde X_t(x)|\big)
  \leq \int_0^T\frac{\big|b_t(X_t(x))-\tilde b_t(\tilde X_t(x))\big|}
  {\rho\big(|X_t(x)-\tilde X_t(x)|\big)+\delta}\,\d t.$$
Since the function $\psi_\delta$ is continuous, we arrive at
  $$\psi_\delta\big(\|X_\cdot(x)-\tilde X_\cdot(x)\|_{\infty,T}\big)
  \leq \int_0^T\frac{\big|b_t(X_t(x))-\tilde b_t(\tilde X_t(x))\big|}
  {\rho\big(|X_t(x)-\tilde X_t(x)|\big)+\delta}\,\d t.$$
Therefore
  \begin{equation}\label{sect-3-thm-1.0}
  \int_{B(R)}\psi_\delta\big(\|X_\cdot(x)-\tilde
  X_\cdot(x)\|_{\infty,T}\big)\,\d x
  \leq \int_0^T\!\!\int_{B(R)}\frac{\big|b_t(X_t(x))-\tilde b_t(\tilde X_t(x))\big|}
  {\rho\big(|X_t(x)-\tilde X_t(x)|\big)+\delta}\,\d x\d t.
  \end{equation}

Denote by $I$ the integral on the right hand side of
\eqref{sect-3-thm-1.0}. Using the triangular inequality, we obtain
  \begin{equation}\label{sect-3-thm-1.1}
  I\leq \int_0^T\!\!\int_{B(R)}\frac{\big|b_t(X_t(x))-b_t(\tilde X_t(x))\big|}
  {\rho\big(|X_t(x)-\tilde X_t(x)|\big)+\delta}\,\d x\d t
  +\int_0^T\!\!\int_{B(R)}\frac{\big|b_t(\tilde X_t(x))-\tilde b_t(\tilde X_t(x))\big|}
  {\rho\big(|X_t(x)-\tilde X_t(x)|\big)+\delta}\,\d x\d t.
  \end{equation}
Since the flows $X_t$ and $\tilde X_t$ leave the Lebesgue measure
absolutely continuous, we can apply the condition (H) for $b_t$ and
obtain that for a.e. $x\in B(R)$,
  $$\big|b_t(X_t(x))-b_t(\tilde X_t(x))\big|\leq \big(g_t(X_t(x))+g_t(\tilde X_t(x))\big)
  \rho\big(|X_t(x)-\tilde X_t(x)|\big).$$
Next it is clear that from Definition \ref{sect-2-def}(i), one has
$\|X_\cdot(x)\|_{\infty,T}\leq R+T\|b\|_{L^\infty}$ and $\|\tilde
X_\cdot(x)\|_{\infty,T}\leq R+T\|\tilde b\|_{L^\infty}$ for a.e.
$x\in B(R)$. Therefore, by the definition of the compressibility
constants $L$ and $\tilde L$, the first term on the right hand side
of \eqref{sect-3-thm-1.1} can be estimated as follows:
  \begin{align}\label{sect-3-thm-1.2}
  \int_0^T\!\!\int_{B(R)}\frac{\big|b_t(X_t(x))-b_t(\tilde X_t(x))\big|}
  {\rho\big(|X_t(x)-\tilde X_t(x)|\big)+\delta}\,\d x\d t
  &\leq\int_0^T\!\!\int_{B(R)}\big(g_t(X_t(x))+g_t(\tilde X_t(x))\big)\,\d x\d t\cr
  &\leq L\int_0^T\!\!\int_{B(\bar R)}g_t(y)\,\d y\d t
  +\tilde L\int_0^T\!\!\int_{B(\bar R)}g_t(y)\,\d y\d t\cr
  &= (L+\tilde L)\|g\|_{L^1([0,T]\times B(\bar R))}.
  \end{align}
Moreover, the second integral in \eqref{sect-3-thm-1.1} is dominated
by
  \begin{align*}
  \frac1\delta\int_0^T\!\!\int_{B(R)}\big|b_t(\tilde X_t(x))-\tilde b_t(\tilde
  X_t(x))\big|\,\d x\d t
  &\leq \frac{\tilde L}\delta \int_0^T\!\!\int_{B(\bar R)}
  |b_t(y)-\tilde b_t(y)|\,\d y\d t\cr
  &= \frac{\tilde L}\delta\|b-\tilde b\|_{L^1([0,T]\times B(\bar
  R))}.
  \end{align*}
Combining this with \eqref{sect-3-thm-1.1} and
\eqref{sect-3-thm-1.2}, we obtain
  $$I\leq (L+\tilde L)\|g\|_{L^1([0,T]\times B(\bar R))}
  +\frac{\tilde L}\delta\|b-\tilde b\|_{L^1([0,T]\times B(\bar
  R))}.$$
Substituting $I$ into \eqref{sect-3-thm-1.0}, we complete the proof.
\fin

\medskip

Now we can prove the main result of this section.

\begin{theorem}[Existence and uniqueness of the regular Lagrangian
flow]\label{sect-3-thm-2}

Assume that $b:[0,T]\times\R^d\ra\R^d$ is a bounded vector field
satisfying {\rm(H)} with $g\in L^1\big([0,T],L_{loc}^1(\R^d)\big)$.
Moreover, the distributional divergence $\div(b)$ of $b$ exists and
$[\div(b)]^- \in L^1\big([0,T],L^\infty(\R^d)\big)$. Then there
exists a unique regular Lagrangian flow generated by $b$.
\end{theorem}

\noindent{\bf Proof.} {\it Step 1: Uniqueness.} Suppose there are
two regular Lagrangian flows $X_t$ and $\hat X_t$ associated to $b$
with compressibility constants $L$ and $\hat L$ respectively.
Applying Theorem \ref{a-priori estimate}, we have
  \begin{equation}\label{sect-3-thm-2.0}
  \int_{B(R)}\psi_\delta\big(\|X_\cdot(x)-\hat X_\cdot(x)\|_{\infty,T}\big)\,\d x
  \leq (L+\hat L)\|g\|_{L^1([0,T]\times B(\bar R))},
  \end{equation}
where $\bar R=R+T\|b\|_{L^\infty}$. If $\L^d\{x\in
B(R):X_\cdot(x)\neq \hat X_\cdot(x)\}>0$, then there is $\ee_0>0$
such that $K_{\ee_0}:=\{ x\in B(R):\|X_\cdot(x)-\hat
X_\cdot(x)\|_{\infty,T}\geq \ee_0\}$ has positive measure. Thus by
\eqref{sect-3-thm-2.0},
  $$(L+\hat L)\|g\|_{L^1([0,T]\times B(\bar R))}\geq
  \int_{K_{\ee_0}}\psi_\delta\big(\|X_\cdot(x)-\hat X_\cdot(x)\|_{\infty,T}\big)\,\d x
  \geq \psi_\delta(\ee_0)\L^d(K_{\ee_0}).$$
Letting $\delta\da0$, we deduce from Lemma \ref{sect-2-lem.1}(1)
that
  $$\infty>(L+\hat L)\|g\|_{L^1([0,T]\times B(\bar R))}\geq \infty,$$
which is a contradiction. Hence $N=\{x\in B(R):X_\cdot(x)\neq \hat
X_\cdot(x)\}$ is $\L^d$-negligible. We conclude that the two flows
$X_t$ and $\hat X_t$ coincide with each other on the interval
$[0,T]$.

{\it Step 2: Existence.} Let $\{\chi_n:n\geq1\}$ be a sequence of
standard convolution kernels. For $t\in[0,T]$, define
$b^n_t=b_t\ast\chi_n$, i.e. the convolution of $b_t$ and $\chi_n$.
Then for every $n\geq1$, $b^n$ is a time dependent smooth vector
field, and
  $$\|b^n_t\|_{L^\infty}\leq \|b_t\|_{L^\infty},\quad
  \big\|[\div(b^n_t)]^-\big\|_{L^\infty} \leq \big\|[\div(b_t)]^-\big\|_{L^\infty},\quad t\in [0,T].$$
Let $X^n_t$ be the smooth flow generated by $b^n_t$, then it is easy
to know that $(X^n_t)_\#\L^d\leq L_n\L^d$, where
  $$L_n=\exp\bigg(\int_0^T\big\|[\div(b^n_t)]^-\big\|_{L^\infty}\,\d t\bigg)
  \leq \exp\bigg(\int_0^T\big\|[\div(b_t)]^-\big\|_{L^\infty}\,\d t\bigg)=:L.$$

Now we show that $b^n_t$ satisfies (H) with $g^n_t=g_t\ast\chi_n$
for all $n\geq1$. To this end, we fix any two points $x,y\in \R^d$.
We have by the definition of convolution,
  $$|b^n_t(x)-b^n_t(y)|\leq \int_{\R^d}|b_t(x-z)-b_t(y-z)|\chi_n(z)\,\d z.$$
Now we shall make use of the condition (H). Note that
$(x-N_t)\cup(y-N_t)$ is a negligible subset. When $z\notin
(x-N_t)\cup(y-N_t)$, one has $x-z\notin N_t$ and $y-z\notin N_t$,
hence by (H),
  $$|b_t(x-z)-b_t(y-z)|\leq \big(g_t(x-z)+g_t(y-z)\big)\rho(|x-y|).$$
As a result,
  \begin{align}\label{sect-3-thm-2.1}
  |b^n_t(x)-b^n_t(y)|&\leq\int_{\R^d}\big(g_t(x-z)+g_t(y-z)\big)
  \rho(|x-y|)\chi_n(z)\,\d z\cr
  &= \big(g^n_t(x)+g^n_t(y)\big)\rho(|x-y|).
  \end{align}
Thus $b^n_t$ satisfies (H) with the function $g^n_t$. Notice that
\eqref{sect-3-thm-2.1} holds for all $x,y\in\R^d$.

From the above arguments, we can apply Theorem \ref{a-priori
estimate} to the sequence of smooth flows $\{X^n_t:n\geq1\}$ and get
  \begin{align}\label{sect-3-thm-2.2}
  &\int_{B(R)}\psi_\delta\big(\|X^n_\cdot(x)-X^m_\cdot(x)\|_{\infty,T}\big)\,\d x\cr
  &\hskip6mm\leq (L_n+L_m)\|g_n\|_{L^1([0,T]\times B(\bar R))}
  +\frac{L_m}\delta \|b^n-b^m\|_{L^1([0,T]\times B(\bar R))}\cr
  &\hskip6mm\leq 2L\|g\|_{L^1([0,T]\times B(\bar R+1))}
  +\frac{L}\delta \|b^n-b^m\|_{L^1([0,T]\times B(\bar R))}.
  \end{align}
Set
  $$\delta=\delta_{n,m}:=\|b^n-b^m\|_{L^1([0,T]\times B(\bar R))}$$
which tends to 0 as $n,m\ra\infty$, since $b^n$ converges to $b$ in
$L^1\big([0,T],L^1_{loc}(\R^d)\big)$. Then we obtain
  \begin{equation}\label{sect-3-thm-2.3}
  \int_{B(R)}\psi_{\delta_{n,m}}\big(\|X^n_\cdot(x)-X^m_\cdot(x)\|_{\infty,T}\big)\,\d x
  \leq 2L\|g\|_{L^1([0,T]\times B(\bar R+1))}+L=:C<\infty.
  \end{equation}

We will show that $\{X^n_\cdot: n\geq1\}$ is a Cauchy sequence in
$L^1(B(R),\Gamma_T)$. For any $\eta>0$, let
  \begin{align*}
  K_{n,m}&=\{x\in B(R):\|X^n_\cdot(x)-X^m_\cdot(x)\|_{\infty,T}\leq \eta\}\cr
  &=\big\{x\in B(R):\psi_{\delta_{n,m}}\big(\|X^n_\cdot(x)-X^m_\cdot(x)\|_{\infty,T}\big)
  \leq \psi_{\delta_{n,m}}(\eta)\big\}.
  \end{align*}
By Chebyshev's inequality and \eqref{sect-3-thm-2.3},
  \begin{equation*}
  \L^d(B(R)\setminus K_{n,m})\leq \frac1{\psi_{\delta_{n,m}}(\eta)}
  \int_{B(R)}\psi_{\delta_{n,m}}\big(\|X^n_\cdot(x)-X^m_\cdot(x)\|_{\infty,T}\big)\,\d x
  \leq \frac C{\psi_{\delta_{n,m}}(\eta)}.
  \end{equation*}
Therefore
  \begin{align*}
  \int_{B(R)}\|X^n_\cdot(x)-X^m_\cdot(x)\|_{\infty,T}\,\d x
  &=\bigg(\int_{K_{n,m}}+\int_{B(R)\setminus K_{n,m}}\bigg)
  \|X^n_\cdot(x)-X^m_\cdot(x)\|_{\infty,T}\,\d x\cr
  &\leq \eta \L^d(K_{n,m})+2(R+T\|b\|_{L^\infty})\L^d(B(R)\setminus
  K_{n,m})\cr
  &\leq \eta
  \L^d(B(R))+2\bar R \frac{C}{\psi_{\delta_{n,m}}(\eta)}.
  \end{align*}
Note that as $n,m\ra\infty$, $\delta_{n,m}\ra0$, thus by Lemma
\ref{sect-2-lem.1}(1), $\psi_{\delta_{n,m}}(\eta)\ra\infty$ for any
$\eta>0$. Consequently,
  $$\varlimsup_{n,m\ra\infty}\int_{B(R)}\|X^n_\cdot(x)-X^m_\cdot(x)\|_{\infty,T}\,\d x
  \leq \eta \L^d(B(R)).$$
By the arbitrariness of $\eta>0$, we conclude that $\{X^n_\cdot:
n\geq1\}$ is a Cauchy sequence in $L^1(B(R),\Gamma_T)$ for any
$R>0$. Therefore, there exists a measurable map $X:\R^d \ra\Gamma_T$
which is the limit in $L^1_{loc}(\R^d,\Gamma_T)$ of $X^n$. We can
find a subsequence $\{n_k:k\geq1\}$, such that for a.e. $x\in\R^d$,
$X^{n_k}_t(x)$ converges to $X_t(x)$ uniformly in $t\in [0,T]$.
Hence we still have
  \begin{equation}\label{sect-3-thm-2.4}
  \|X_\cdot(x)\|_{\infty,T}\leq R+T\|b\|_{L^\infty}=\bar R,\quad\mbox{for a.e. }x\in B(R).
  \end{equation}

Now we prove that $X_t$ is a regular Lagrangian flow generated by
$b$. Firstly, for any $\phi\in C_c(\R^d,\R_+)$, we have by the Fatou
lemma,
  $$\int_{\R^d}\phi(X_t(x))\,\d x\leq \varliminf_{k\ra\infty}
  \int_{\R^d}\phi(X^{n_k}_t(x))\,\d x\leq L\int_{\R^d}\phi(y)\,\d y.$$
This implies
  \begin{equation}\label{sect-3-thm-2.5}
  (X_t)_\#\L^d\leq L\L^d,\quad \mbox{for all }t\in[0,T];
  \end{equation}
thus Definition \ref{sect-2-def}(ii) is satisfied. Secondly, we show
that for $\L^d$-a.e. $x\in\R^d$, $t\ra X_t(x)$ is an integral curve
of the vector field $b_t$. To this end, we estimate the quantity
  $$J^n:=\int_{B(R)}\sup_{0\leq t\leq T}\bigg|\int_0^tb^n_s(X^n_s(x))\,\d s
  -\int_0^tb_s(X_s(x))\,\d s\bigg|\,\d x.$$
By the triangular inequality, $J^n$ is dominated by the sum of
  $$J^n_1:=\int_{B(R)}\!\int_0^T\big|b^n_s(X^n_s(x))-b_s(X^n_s(x))\big|\,\d s\d x$$
and
  $$J^n_2:=\int_{B(R)}\!\int_0^T\big|b_s(X^n_s(x))-b_s(X_s(x))\big|\,\d s\d x.$$
For the first term, we have
  $$J^n_1=\int_0^T\!\!\int_{B(R)}\big|b^n_s(X^n_s(x))-b_s(X^n_s(x))\big|\,\d x\d s
  \leq L\int_0^T\!\!\int_{B(\bar R)}\big|b^n_s(y)-b_s(y)\big|\,\d y\d s.$$
Hence
  \begin{equation}\label{sect-3-thm-2.6}
  \lim_{n\ra\infty}J^n_1=0.
  \end{equation}
For any $\ee>0$, we take a vector field $\hat b\in C^1([0,T]\times
B(\bar R),\R^d)$ such that
  $$\int_0^T\!\!\int_{B(\bar R)}\big|\hat b_s(x)-b_s(x)\big|\,\d x\d s<\ee.$$
Again by the triangular inequality,
  \begin{align*}
  J^n_2&\leq \int_0^T\!\!\int_{B(R)}\big|b_s(X^n_s(x))-\hat b_s(X^n_s(x))\big|\,\d x\d s
  +\int_0^T\!\!\int_{B(R)}\big|\hat b_s(X^n_s(x))-\hat b_s(X_s(x))\big|\,\d x\d
  s\cr
  &\hskip11pt +\int_0^T\!\!\int_{B(R)}\big|\hat b_s(X_s(x))-b_s(X_s(x))\big|\,\d x\d
  s\cr
  &=:J^n_{2,1}+J^n_{2,2}+J^n_{2,3}.
  \end{align*}
Since $(X^n_s)_\#\L^d\leq L\L^d$ for all $s\in[0,T]$, we have
  $$J^n_{2,1}\leq L\int_0^T\!\!\int_{B(\bar R)}\big|b_s(y)-\hat b_s(y)\big|\,\d y\d s
  <L\ee.$$
By \eqref{sect-3-thm-2.4} and \eqref{sect-3-thm-2.5}, the same
argument leads to
  $$J^n_{2,3}<L\ee.$$
Moreover, by the choice of $\hat b$, there is $C_1>0$ such that
$\sup_{0\leq s\leq T}\big\|\nabla\hat b_s\big\|_{L^\infty(B(\bar
R))}\leq C_1$. Therefore
  $$J^n_{2,2}\leq C_1\int_0^T\!\!\int_{B(R)}|X^n_s(x)-X_s(x)|\,\d x\d s
  \leq C_1T\int_{B(R)}\|X^n_\cdot(x)-X_\cdot(x)\|_{\infty,T}\,\d x\ra0$$
as $n$ goes to $\infty$. Summing up the above arguments, we get
  $$\varlimsup_{n\ra\infty}J^n_2\leq 2L\ee.$$
Since $\ee>0$ is arbitrary, we obtain $\lim_{n\ra\infty}J^n_2=0$.
Combining this with \eqref{sect-3-thm-2.6}, we finally obtain
$\lim_{n\ra\infty}J^n=0$. Therefore, letting $n\ra\infty$ in the
equality
  $$X^n_t(x)=x+\int_0^t b^n_s(X^n_s(x))\,\d s\quad\mbox{for all }t\leq T,$$
we see that both sides converge in $L^1_{loc}(\R^d,\Gamma_T)$ to
$X_\cdot(x)$ and $x+\int_0^\cdot b_s(X_s(x))\,\d s$ respectively.
Hence for $\L^d$-a.e. $x\in\R^d$, it holds
  $$X_t(x)=x+\int_0^t b_s(X_s(x))\,\d s\quad\mbox{for all }t\in[0,T];$$
that is, $t\ra X_t(x)$ is an integral curve of the vector field
$b_t$. To sum up, $X_t$ is a regular Lagrangian flow generated by
$b$. \fin

\begin{remark} The condition $[\div(b)]^- \in L^1\big([0,T],
L^\infty(\R^d)\big)$ can be relaxed as $[\div(b)]^- \in
L^1\big([0,T],L^\infty_{loc}(\R^d)\big)$, since we have good
estimate on the growth of the flow $X_t$.
\end{remark}

\begin{remark} Under the condition (H), it seems to the authors that one is unable
to prove the well posedness of the transport equation
  \begin{equation*}
  \frac{\partial}{\partial t}u_t+b_t\cdot\nabla u_t +c_t\,u_t=0,\quad
  u|_{t=0}=u_0.
  \end{equation*}
by following DiPerna-Lions's original approach, that is, by showing
that the commutator
  $$r_n(b_t,u_t)=(b_t\cdot\nabla u_t)\ast\chi_n-b_t\cdot\nabla(u_t\ast\chi_n)$$
converges to 0 strongly in $L^1_{loc}$. Here $\chi_n$ is the
standard convolution kernel. This can be seen from the proof of
\cite[Lemma II.1]{DiPernaLions89} (or \cite[Proposition
4.1]{Ambrosio08}), which essentially relies on the ``almost
everywhere Lipschitz continuity'' of Sobolev vector fields.
\end{remark}

\section{Regularity of the flow}

In this section, we first prove a regularity result on the regular
Lagrangian flow, a property much weaker than the approximate
differentiability discussed in \cite{CrippadeLellis}. We need the
following notation: for a bounded measurable subset $U$ with
positive measure, define the average of $f\in L^1_{loc}(\R^d)$ on
$U$ by
  $$\bint_{U}f\,\d x=\frac1{\L^d(U)}\int_U f\,\d x.$$
Then the local maximal function
  $$M_Rf(x)=\sup_{0<r\leq R}\bint_{B(x,r)}|f(y)|\,\d y.$$

Now we can prove

\begin{theorem}\label{sect-4-thm-1}
Let $b$ be a bounded vector field satisfying {\rm(H)}, and
$[\div(b)]^- \in L^1\big([0,T], L^\infty(\R^d)\big)$. Let $X$ be the
unique regular Lagrangian flow associated to $b$. Then for any $R>0$
and sufficiently small $\ee$, there are a measurable subset
$E\subset B(R)$ and some constant $C$ depending on $R, d$ and $g$,
such that $\L^d(B(R)\setminus E)\leq\ee$ and for all $t\in[0,T]$ and
$x,y\in E$, one has
  $$|X_t(x)-X_t(y)|\leq \psi_{|x-y|}^{-1}(C/\ee).$$
\end{theorem}

Here $\psi_r^{-1}$ is the inverse function of $\psi_r$. Note that by
Lemma \ref{sect-2-lem.1}(1), we have
$\lim_{r\da0}\psi^{-1}_r(\xi)=0$ for all $\xi>0$. Therefore this
theorem implies that $X_t$ is uniformly continuous in $E$, since
when $y\ra x$ in the subset $E$, the quantity
$\psi_{|x-y|}^{-1}(C/\ee)$ decreases to 0. Unfortunately, the
function $\psi_r^{-1}$ does not have an explicit expression, unless
$\rho(s)=s$ for all $s\geq0$ (see Remark \ref{sect-4-rem}).

\medskip

\noindent{\bf Proof of Theorem \ref{sect-4-thm-1}.} We follow the
ideas of \cite[Remark 2.4]{CrippadeLellis} (see also
\cite[Proposition 5.2]{LiLuo}). For $0\leq t\leq T$, $0<r\leq 2R$
and $x\in B(R)$, define
  $$Q(t,x,r)=\bint_{B(x,r)}\psi_r(|X_t(x)-X_t(y)|)\,\d y.$$
Then
  $$Q(0,x,r)=\bint_{B(x,r)}\psi_r(|x-y|)\,\d y
  \leq \bint_{B(x,r)}\psi_r(r)\,\d y\leq 1.$$
By Definition \ref{sect-2-def}(i), we see that $t\mapsto Q(t, x,r)$ is
Lipschitz and
  \begin{align*}
  \frac{\d}{\d t}Q(t,x,r)&=\bint_{B(x,r)}\psi'_r(|X_t(x)-X_t(y)|)
  \,\frac{\d}{\d t}|X_t(x)-X_t(y)|\,\d y\cr
  &\leq\bint_{B(x,r)}\frac{\big|b_t(X_t(x))-b_t(X_t(y))\big|}
  {\rho(|X_t(x)-X_t(y)|)+r}\,\d y.
  \end{align*}
Using the assumption (H) on $b$, we have
  $$\frac{\d}{\d t}Q(t,x,r)\leq \bint_{B(x,r)}\big(g_t(X_t(x))+g_t(X_t(y))\big)\,\d y
  =g_t(X_t(x))+\bint_{B(x,r)}g_t(X_t(y))\,\d y.$$
Integrating both sides with respect to time from $0$ to $t$, we
arrive at
  \begin{align}\label{sect-4-thm-1.0}
  Q(t,x,r)&\leq Q(0,x,r)+\int_0^t g_s(X_s(x))\,\d s
  +\int_0^t\!\!\bint_{B(x,r)}g_s(X_s(y))\,\d y\d s\cr
  &\leq 1+\int_0^T g_s(X_s(x))\,\d s
  +\int_0^T\!\!\bint_{B(x,r)}g_s(X_s(y))\,\d y\d s.
  \end{align}
Denote by $\Phi(x)=\int_0^T g_s(X_s(x))\,\d s$ for a.e. $x\in\R^d$.
Then for all $t\leq T$,
  $$Q(t,x,r)\leq 1+\Phi(x)+\bint_{B(x,r)}\Phi(y)\,\d y.$$
Therefore
  \begin{equation}\label{sect-4-thm-1.1}
  \sup_{0\leq t\leq T}\sup_{0<r\leq 2R}Q(t,x,r)\leq 1+\Phi(x)+M_{2R}\Phi(x).
  \end{equation}

For $\eta>0$ sufficiently small, we have
  \begin{align*}
  I&:=\L^d\{x\in B(R):1+\Phi(x)+M_{2R}\Phi(x)>1/\eta\}\cr
  &\leq \L^d\{x\in B(R):\Phi(x)>1/(3\eta)\}
  +\L^d\{x\in B(R):M_{2R}\Phi(x)>1/(3\eta)\}.
  \end{align*}
By Chebyshev's inequality and Lemma \ref{sect-2-lem.2}, we have
  \begin{align*}
  I&\leq 3\eta\int_{B(R)}\Phi(x)\,\d x+3\eta
  C_d\int_{B(3R)}\Phi(y)\,\d y\cr
  &\leq 3\eta(1+C_d)\int_{B(3R)}\Phi(y)\,\d y.
  \end{align*}
By the definition of $\Phi$, one has
  \begin{align*}
  I&\leq 3\eta(1+C_d)\int_0^T\!\!\int_{B(3R)}g_t(X_t(y))\,\d y\d
  t\cr
  &\leq 3\eta(1+C_d)L\int_0^T\!\!\int_{B(3R+T\|b\|_{L^\infty})}
  g_t(x)\,\d x\d t.
  \end{align*}
Let $\bar C:=3(1+C_d)L\|g\|_{L^1([0,T]\times
B(3R+T\|b\|_{L^\infty}))}$; then $I\leq \eta \bar C$.

Now for any $\ee>0$, set $\eta=\ee/\bar C$. Then by
\eqref{sect-4-thm-1.1} and the definition of $I$,\
  $$\L^d\bigg\{x\in B(R):\sup_{0\leq t\leq T}\sup_{0<r\leq 2R}Q(t,x,r)
  >\frac{\bar C}\ee\bigg\}\leq I\leq \eta \bar C=\ee.$$
Let
  $$E=\bigg\{x\in B(R):\sup_{0\leq t\leq T}\sup_{0<r\leq 2R}Q(t,x,r)
  \leq\frac{\bar C}\ee\bigg\}.$$
Then $\L^d(B(R)\setminus E)\leq \ee$ and for any $x\in E$, $0\leq
t\leq T$ and $0<r\leq 2R$, one has
  \begin{equation}\label{sect-4-thm-1.2}
  \bint_{B(x,r)}\psi_r(|X_t(x)-X_t(y)|)\,\d y\leq\frac{\bar C}\ee.
  \end{equation}

Now fix any $x,y\in E$ and let $r=|x-y|$ which is less than $2R$.
Lemma \ref{sect-2-lem.1}(2) tells us that $\psi_r$ is concave, hence
$\psi_r(a+b)\leq\psi_r(a)+\psi_r(b)$ for any $a,b\geq0$. As a
result,
  $$\psi_r(|X_t(x)-X_t(y)|)\leq \psi_r(|X_t(x)-X_t(z)|)+\psi_r(|X_t(z)-X_t(y)|).$$
Therefore,
  \begin{align*}
  \psi_r(|X_t(x)-X_t(y)|)
  &=\bint_{B(x,r)\cap B(y,r)}\psi_r(|X_t(x)-X_t(y)|)\,\d z\cr
  &\leq\bint_{B(x,r)\cap B(y,r)}\psi_r(|X_t(x)-X_t(z)|)\,\d z\cr
  &\hskip11pt+\bint_{B(x,r)\cap B(y,r)}\psi_r(|X_t(z)-X_t(y)|)\,\d
  z.
  \end{align*}
Let $\tilde C_d=\L^d(B(x,r))/\L^d(B(x,r)\cap B(y,r))$ which only
depends on the dimension $d$; then
  \begin{align*}
  \psi_r(|X_t(x)-X_t(y)|)
  &\leq \tilde C_d\bint_{B(x,r)}\psi_r(|X_t(x)-X_t(z)|)\,\d z
  +\tilde C_d\bint_{B(y,r)}\psi_r(|X_t(z)-X_t(y)|)\,\d z\cr
  &\leq 2\tilde C_d\bar C/\ee,
  \end{align*}
where the last inequality follows from \eqref{sect-4-thm-1.2}.
Consequently, for all $t\leq T$ and $x,y\in E$,
  $$|X_t(x)-X_t(y)|\leq\psi_r^{-1}(2\tilde C_d\bar C/\ee)
  =\psi_{|x-y|}^{-1}\big(2\tilde C_d\bar C/\ee\big).$$\fin

\begin{remark}\label{sect-4-rem}
If $\rho(s)=s$ for all $s\geq0$, then $\psi_r(s)=\log \big(\frac sr
+1\big)$ and $\psi_r^{-1}(t)=r(e^t-1)$. Thus the last inequality in
the proof of Theorem \ref{sect-4-thm-1} becomes
  $$|X_t(x)-X_t(y)|\leq |x-y|\big(e^{2C_d\bar C/\ee}-1\big)
  \leq |x-y| e^{2C_d\bar C/\ee},$$
which is the estimate given in \cite[Proposition
2.3]{CrippadeLellis} in the case $p=1$.
\end{remark}

We complete this section by discussing the compactness of the
regular Lagrangian flow, following the ideas in \cite[Section
4]{CrippadeLellis}. For fixed $R>0$ and $0<r<R/2$, set
  $$a(r,R,X)=\int_{B(R)}\sup_{0\leq t\leq T}\bint_{B(x,r)}
  \psi_r(|X_t(x)-X_t(y)|)\,\d y\d x.$$

\begin{proposition}\label{sect-4-prop}
Let $b$ be a bounded vector field satisfying the condition {\rm (H)}
with some function $g\in L^1([0,T],L^1_{loc}(\R^d))$. Let $X$ be a
regular Lagrangian flow associated to $b$ with the compressibility
constant $L$. Then
  $$a(r,R,X)\leq \L^d(B(R))+2L\|g\|_{L^1([0,T]\times B(\bar R))},$$
where $\bar R=3R/2+2T\|b\|_{L^\infty}$.
\end{proposition}

\noindent{\bf Proof.} Recall the definition of $Q(t,x,r)$ at the
beginning of the proof of Theorem \ref{sect-4-thm-1}. We still have
the estimate \eqref{sect-4-thm-1.0}: for any $t\leq T$,
  \begin{align*}
  Q(t,x,r)&\leq 1+\int_0^T g_s(X_s(x))\,\d s
  +\int_0^T\!\! \bint_{B(x,r)}g_s(X_s(y))\,\d y\d s.
  \end{align*}
Therefore
  \begin{align*}
  a(r,R,X)&=\int_{B(R)}\sup_{0\leq t\leq T}Q(t,x,r)\,\d x\cr
  &\leq \L^d(B(R))+\int_{B(R)}\!\int_0^T g_s(X_s(x))\,\d s\d x
  +\int_{B(R)}\!\int_0^T\!\!\bint_{B(x,r)}g_s(X_s(y))\,\d y\d s\d x.
  \end{align*}
We denote by $I_1$ and $I_2$ the two integrals on the right hand
side respectively. First
  $$I_1=\int_0^T\!\!\int_{B(R)} g_s(X_s(x))\,\d x\d s
  \leq \int_0^T L \int_{B(\bar R)} g_s(y)\,\d y\d s
  =L\|g\|_{L^1([0,T]\times B(\bar R))}.$$
For the second integral, by changing the order of integration, we
have
  \begin{align*}
  I_2&=\int_{B(R)}\!\int_0^T\!\!\bint_{B(r)}g_s(X_s(x+z))\,\d z\d s\d x
  =\bint_{B(r)}\!\int_0^T\!\!\int_{B(R)}g_s(X_s(x+z))\,\d x\d s\d z.
  \end{align*}
Therefore,
  $$I_2\leq \bint_{B(r)}\!\int_0^T L\int_{B(\bar R)}g_s(y)\,\d y\d s\d z
  =L\|g\|_{L^1([0,T]\times B(\bar R))}.$$
Combining the above two estimates, we arrive at the conclusion. \fin

\medskip

Now we can prove

\begin{theorem}[Compactness of the flow]\label{sect-4-thm-2}
Let $\{b^n:n\geq1\}$ be a sequence of vector fields equi-bounded in
$L^\infty([0,T]\times\R^d)$. For every $n\geq1$, assume that $b^n$
satisfies {\rm (H)} with the function $g^n$, and the family
$\{g^n:n\geq1\}$ is equi-bounded in $L^1([0,T],L^1_{loc}(\R^d))$.
Let $X^n$ be a regular Lagrangian flow associated to $b^n$ with the
compressibility constant $L_n$. Suppose $\sup_{n\geq1}L_n\leq
L<\infty$. Then the sequence $\{X^n:n\geq1\}$ is strongly precompact
in $L^1_{loc}([0,T]\times\R^d)$.
\end{theorem}

\noindent{\bf Proof.} Applying the estimate in Proposition
\ref{sect-4-prop} to the flow $X^n$, we get
  $$a(r,R,X^n)\leq \L^d(B(R))+2L_n\|g^n\|_{L^1([0,T]\times B(\bar R_n))},$$
where $\bar R_n=3R/2+2T\|b^n\|_{L^\infty}$. Since $\{b^n\}$ is
equi-bounded, we see that $\tilde R:=\sup_{n\geq1}\bar R_n <\infty$.
Moreover, by the boundedness of the sequence $\{g^n\}$ in
$L^1([0,T],L^1_{loc}(\R^d))$, we obtain
  \begin{equation}\label{sect-4-thm-2.1}
  \sup_{n\geq1} a(r,R,X^n)\leq \L^d(B(R))
  +2L\sup_{n\geq1}\|g^n\|_{L^1([0,T]\times B(\tilde R))}=:C_{d,R,T}<\infty.
  \end{equation}
For $0< z\leq \tilde R$, by Lemma \ref{sect-2-lem.1}(2), one has
  $$\frac{\psi_r(z)}z\geq \frac{\psi_r(\tilde R)}{\tilde R},\quad
  \mbox{or equivalently,}\quad z\leq \frac{\tilde R}{\psi_r(\tilde R)}\,\psi_r(z).$$
Since $y\in B(x,r)$ and $r\leq R/2$, it holds
$|X^n_t(x)-X^n_t(y)|\leq\bar R_n\leq \tilde R$. Hence by
\eqref{sect-4-thm-2.1},
  \begin{align*}
  &\int_{B(R)}\sup_{0\leq t\leq T}\bint_{B(x,r)}
  |X^n_t(x)-X^n_t(y)|\,\d y\d x\cr
  &\hskip6mm \leq \frac{\tilde R}{\psi_r(\tilde R)}
  \int_{B(R)}\sup_{0\leq t\leq T}\bint_{B(x,r)}
  \psi_r(|X^n_t(x)-X^n_t(y)|)\,\d y\d x\cr
  &\hskip6mm =\frac{\tilde R}{\psi_r(\tilde R)}\,a(r,R,X^n)
  \leq \frac{\tilde R}{\psi_r(\tilde R)}\,C_{d,R,T}=:g(r),
  \end{align*}
where the function $g(r)$ does not depend on $n$ and satisfies
$g(r)\da0$ as $r$ decreases to 0, by Lemma \ref{sect-2-lem.1}(1).
Similar to the estimate of $I_2$ in the proof of Proposition
\ref{sect-4-prop}, we change the order of integration and obtain
  \begin{equation}\label{sect-4-thm-2.2}
  \sup_{n\geq1}\sup_{0\leq t\leq T}\int_{B(r)}\!\int_{B(R)}
  |X^n_t(x)-X^n_t(x+z)|\,\d x\d z\leq g(r)\L^d(B(r)).
  \end{equation}
The rest of the proof is similar to that of \cite[Corollary
4.2]{CrippadeLellis}, hence we omit it. \fin

\medskip

With the above compactness result in mind, we can give another proof
of the existence of the regular Lagrangian flow.

\begin{corollary}[Existence of the flow]
Let $b$ be a bounded vector field satisfying {\rm (H)} with the
function $g$. Assume that $[\div(b)]^-\in L^1([0,T],
L^\infty(\R^d))$. Then there exists a regular Lagrangian flow
associated to $b$.
\end{corollary}

\noindent{\bf Proof.} We regularize the vector field $b$ as in {\it
Step 2} of the proof of Theorem \ref{sect-3-thm-2}. It is clear that
the conditions of Theorem \ref{sect-4-thm-2} are satisfied by the
smooth vector fields $\{b^n\}$ and the corresponding flows
$\{X^n\}$. As a result, the sequence $\{X^n\}$ is strongly
precompact in $L^1_{loc}([0,T]\times\R^d)$, and every limit point of
$\{X^n\}$ is a regular Lagrangian flow associated to $b$. \fin

\end{document}